\documentclass[12pt]{article}

\author{A. Assi\thanks{Universit\'e d'Angers, D\'epartement de
Math\'ematiques, 2 bd Lavoisier, 49045 Angers Cedex 01, France},  M. Barile \thanks{Universit\`a degli Studi di Bari,
                    Dipartimento Interuniversitario di Matematica,
                    Via E. Orabona 4, 70125 Bari, Italia.\newline The second author would like to thank the Departement of Mathematics of the University of Angers for hospitality and support during the preparation of the present paper.}}

\title{Effective construction of irreducible curve singularities}
\date{\mbox{}}

\newtheorem{teorema}{Theorem}[section]

\newtheorem{proposicion}[teorema]{Proposition}
\newtheorem{lema}[teorema]{Lemma}
\newtheorem{definicion}[teorema]{Definition}

\newtheorem{corolario}[teorema]{Corollary}

\newtheorem{nota}[teorema]{Remark}
\newtheorem{exemple}[teorema]{Example}
\newtheorem{theorema}[teorema]{Theorem}
\newenvironment{demostracion}[1]{\paragraph{\sl Proof#1}}{}

\newcommand{\KK}{{\bf K}}

\newcommand{\ZZ}{{\bf Z}}
\newcommand{\NN}{{\bf N}}

\begin{document}
\maketitle








\begin{center}
 {\bf Abstract} \footnote{ Mathematical Subject Classification:
 14Q05, 14C17, 32S15\newline
Keywords: Curve singularities, approximate roots, Milnor number.}
\end{center}
\noindent {\small
By using the effective notion of the approximate roots of a
polynomial, we describe the equisingularity classes of irreducible
curve singularities with a given Milnor number.}

\section*{Introduction} Let ${\bf K}$ be an algebraically closed field of characteristic 
zero. Let $f$ be an irreducible monic polynomial of
 ${\bf R}={\bf K}[[x]][y]$, say  $f=f(x,y)=y^n+a_1(x).y^{n-1}+\ldots+a_n(x)\in {\bf R}$. Up to a change of coordinates,  we assume that $a_1(x)=0$. 
For all  $ g\in {\bf R}$ let int$(f,g)$ denote the intersection multiplicity of $f$ and $g$. 
Let  $\Gamma(f)=\lbrace {\rm int}(f,g); g\in {\bf R}-(f)\rbrace $ be the semigroup of $f$.  If $f'$ is  another irreducible polynomial of
 ${\bf R}$, then $f$ and $f'$ are said to be 
{\it equisingular} if $\Gamma(f)=\Gamma(f')$ (for example, $y^2-x^3$
and $y^3-x^2$ are equisingular because they are both associated with the  semigroup
generated by $2,3$. In particular, two equisingular polynomials of
${\bf R}$ need not have the same degree in $y$). It is well-known that in this case $\mu(f)=\mu(f')$, where $\mu(f)=int(f_x,f_y)$ is called the
 {\it Milnor number} of $f$.  The converse is  false.   The  {\it
equisingularity class} of the polynomial $f$  is the set of
irreducible polynomials of ${\bf R}$ which are equisingular to $f$. It
is of a certain interest to determine this  equisingularity class,
which gives a classification of the polynomials of ${\bf R}$  in terms
of subsemigroups of $\ZZ$.  Another remarkable classification is
obtained if one  can characterize  all polynomials whose Milnor
number is equal to some  fixed  nonnegative integer $m$. The aim of
this paper is to study the two questions from an effective point of
view: we first give, for a fixed semigroup of an irreducible
polynomial $f$ of ${\bf R}$, all  elements of the equisingularity
class of $f$. Then, for a fixed $m$ in $\NN$,  by similar methods we
construct  the generic forms of all  irreducible polynomials $f$ of
${\bf R}$ such that $\mu(f)=m$. The set of these polynomials is the
union of a finite number of equisingularity classes. We think that this effective 
classification is useful in the study of problems and conjectures 
in the theory of irreducible curve singularities, particularly in the
understand  of their moduli spaces. Our approach 
uses the effective notion of approximate roots of a 
polynomial $f$ of ${\bf R}$ introduced by S.S. Abhyankar and T.T. Moh and 
the notion of generalized Newton polygon introduced by Abhyankar. The first 
one gives rise to an algorithm for the computation of  the set of generators 
of the semigroup of $f$ (and then the set of Newton-Puiseux 
pairs of $f$, see definition 1.3., and [1], [2]). The second one is used by 
Abhyankar to give an irreducibility criterion for the polynomial $f$ (see [3]).

\noindent We would like to point out that our algorithms are intrinsic and that they
have been implemented with {\it Mathematica} (see [8]), and {\it Maple}.

\section{Characteristic sequences}  In this Section we recall the notion of approximate roots of $f$ as well as the characteristic sequences associated with an irreducible polynomial $f=y^n+a_2(x).y^{n-2}+\ldots+a_n(x)$ of ${\bf R}$. 

\begin{definicion} {\rm For any monic polynomial $g\in{\bf R}$,  the {\bf  intersection multiplicity} int$(f,g)$ of $f$ with $g$  is the $x$-order of the $y$-resultant of $f$ and $g$.}\end{definicion}
\noindent
The set $\Gamma(f)=\lbrace {\rm int}(f,g); g\in {\bf R}-f\rbrace $ is a subsemigroup of $\ZZ$, called the {\bf semigroup of $f$}.

\begin{definicion} {\rm Let $y(t)=\sum_j a_jt^j\in
 {\bf K}[[t]]$ be 
 a root of $f(t^n,y)=0$, according to Newton Theorem. Then 
 set $m_0=d_1=n, m_1=$ inf$\lbrace j; a_j\not=0 \rbrace$, and for all $k
\geq 1,$ let

\noindent $m_{k+1}=$ inf$\lbrace j; a_j\not=0$ and $d_k$ does not divide
$j \rbrace$, and $d_{k+1}={\rm gcd}(m_{k+1},d_k)$.

\noindent Since $f$ is irreducible, there exists $h$ such that
$d_{h+1}=1$. We set $m_{h+1}=+\infty$. 

\noindent Finally set $r_0=m_0=n$, $r_1= O_x(a_n(x))$-where $O_x$ 
 denotes the $x$-order-, and for all $k=1,\dots, h-1$:

$r_{k+1}=r_k ({\displaystyle\frac{d_k}{d_{k+1}}}) + (m_{k+1}-m_k)$.}\end{definicion}

\noindent (Remark that, since $a_1(x)=0$, then $r_1=m_1$).

\noindent We recall that with respect to these notations,
$r_0,\dots,r_h$ generates the semigroup $\Gamma(f)$ of $f$. We
denote $\Gamma(f)=<r_0,\ldots,r_h>$.

 \begin{definicion} {\rm For all  $k=1,\ldots, h$, set  
 $e_k={\displaystyle\frac{d_k}{d_{k+1}}}$. The  set   
 $\lbrace ({\displaystyle\frac{m_k}{d_{k+1}}},e_k), 1\leq k\leq h\rbrace$ is called the set 
 of {\bf Newton-Puiseux pairs } of $f$.}\end{definicion}

\begin{definicion} {\rm Let $d$ be a positive integer
and assume that $d$ divides $n$. Let $g$ be a monic polynomial of 
${\bf R}$, of degree ${\displaystyle\frac{n}{d}}$
in $y$. We call $g$ the $d$-th approximate root of $f$ if one of
the following holds:

i) deg$_y (f-g^d) < n-{\displaystyle\frac{n}{d}}$.

ii) in the expansion $f=g^d +\alpha_1 g^{d-1}+\ldots + \alpha_d$ of $f$ with respect to the powers of $g$,
$\alpha_1=0$.

\noindent Remark that i) and ii) are equivalent.

\noindent We denote the $d$-th approximate root of $f$ by
App$_d(f)$. It is clear that App$_d(f)$ is unique, and also that it is
effectively computable if the series $a_k(x), k=2,\ldots,n$, are
polynomials.}
\end{definicion}

\begin{nota}{\rm Given a divisor $d$ of $n$, the 
$d$th approximate root App$_d(f)$ of $f$ can be effectively constructed from
the equation of $f$ in the following way:

\noindent Take $G_0=y^{n/d}$ and let 
$f=G_0^d +\alpha^0_1 G_0^{d-1}+\ldots + \alpha^0_d$ be the expansion
of $f$ with respect to the powers of $G_0$.

i) If $\alpha^0_1=0$, then $G_0= {\rm App}_d(f)$.

ii)  If $\alpha^0_1\not=0$, then set 
$G_1=G_0+\displaystyle{{\alpha^0_1\over d}}$ and consider the expansion           
$f=G_1^d +\alpha^1_1 G_1^{d-1}+\ldots + \alpha^1_d$ of $f$ with respect
to the powers of $G_1$. If $\alpha^1_1\not=0$, then easy calculations show that 
deg$_y\alpha^0_1 > {\rm deg}_y\alpha^1_1$. This procees shall stop 
after a finite number of steps, constructing App$_d(f)$.}
\end{nota}

\begin{nota}{\rm If the characteristic of $\KK$ is not zero and if 
this characteristic does not divide $n$, then the 
construction above applies without any restriction. Otherwise,
the theory of approximate roots does not work as it. Further
information can be found in [9].}
\end{nota}

\noindent Let $g_1,\ldots,g_h, g_{h+1}$ be the  $d_k$-th approximate 
roots of $f$, for $k=1,\ldots,h+1$ (in particular
$g_1=y$ and $g_{h+1}=f$).

\begin{lema} (see [1], (8.2) the Fundamental Theorem (part one)) For all $k=1,\ldots,h$, we have:

i) {\rm int}($f,g_k)=r_k$.


ii) $g_k$ is irreducible in ${\bf R}$ and 
$\Gamma (g_k) = < {\displaystyle\frac{r_0}{d_k}},\ldots,
{\displaystyle\frac{r_{k-1}}{d_k}}>$. Furthermore, 
$g_1,\ldots,$ $g_{k-1}$ are the approximate roots of $g_k$.\end{lema}

\begin{lema} (see [13])  The following formulas hold:

 - {\rm int}$(f_x, f_y)=\displaystyle\sum_{i=1}^{h} (e_i-1)r_i
   -n+1$. In particular ${\rm int}(f_x, f_y)$ is even.

 - For all $k=2,\dots, h$, 
 {\rm int}$(f_x, f_y)=d_k.${\rm int}$(g_{k_x}, g_{k_y})+
 \displaystyle\sum_{i=k}^{h} (e_i-1)r_i - d_k + 1.$
\end{lema}

\begin{demostracion}{.} The proof of the first formula can be found in 
[13] (3.14., page 18). The second formula results from the first one
by easy calculations.
\end{demostracion}

\begin{nota}{\rm The intersection multiplicity int$(f_x,f_y)$ is also called 
the Milnor number of $f$. It is an invariant of $f$ and, by the
formula above, it is common to the elements of the equisingularity
class of $f$. It also coincides with the conductor of the semigroup
$\Gamma(f)$ -usually denoted by $c$- which has the following 
numerical characterization: for all $p\geq c, p\in \Gamma(f)$. Furthermore, given two 
integers $a,b$, if $a+b=m-1$ then exactly one of $a,b\in \Gamma(f)$. It follows that,
since $\Gamma(f)$ has no negative integers, Card(${\bf N}-\Gamma)={\displaystyle\frac{m}{2}}$.
In fact, $c$ is nothing but the order of the conductor of the quotient
$\displaystyle {{\frac {\bf R} {(f)}}}$ into its integral closure. Contrary to
the Milnor number, the conductor can be defined without restriction on the characteristic
of $\KK$. An exhaustive exposition of this theory in positive characteristic
can be found in [9].}
\end{nota}

\bigskip

\section{Generalized Newton polygons and the irreducibility criterion of
Abhyankar}

\medskip

\noindent Let $f= y^n+a_2(x)y^{n-2}+\ldots+a_n(x)$ be a 
monic polynomial, non  necessarily irreducible in ${\bf R}$. In this section 
the notations introduced above will have a more general meaning: 
 $r=(r_0=n,r_1,\ldots,r_h)$ will denote any 
sequence of integers such that $r_k < r_{k+1}$ for all
$k=1,\ldots,h-1$, and we shall set  $d_{k+1}=$ gcd$(r_0,r_1,\ldots, r_k)$ for all
$k=0,\ldots,h$. For all $k=1,\ldots,h$, we set 
$e_k={\displaystyle\frac{d_k}{d_{k+1}}}$;
$g =(g_1,\ldots,g_h,g_{h+1}=f)$ will be a sequence of monic
polynomials of ${\bf R}$ such that deg$_yg_k ={\displaystyle\frac{n}{d_k}}$ for all 
$k=1,\ldots,h$. We recall some important properties.

\medskip

\begin{theorema}{\rm  (see [1], (8.3) The fundamental Theorem (part two)) Let 

$$
B=\lbrace b=(b_1,b_2,\ldots,b_h,b_{h+1})\in \NN^{h+1}; 
b_1<e_1,\ldots,b_h<e_h, b_{h+1} <+\infty\rbrace
$$

\noindent For all $ b\in B$, denote  
$g^{b}=g_1^{b_1}.\ldots.g_h^{b_h}.f^{b_{h+1}}$, then we have:

i) ${\bf R}=\sum_{ b\in B}{\bf K}[[x]]. g^{b}$.

ii) Let $p$ be a polynomial of ${\bf R}$ and 
write $p=\sum_{k=1}^s a_k(x).g^{b^k}$, where $b^k\in B$ for all $k=1,\ldots,s$. Moreover let $b^k_0=$ O$_x a_k(x)$, then 
associate with any ``monomial'' $ a_k(x).g^{b^k}$ the 
integer $< (b^k_0,b^k_1,\ldots,b^k_h), r> = b^k_0.r_0+\sum_{i=1}^{h} b^k_i.r_i$. Finally 
let $B'=\lbrace  b^k; b^k_{h+1}=0\rbrace$. With these notations we have the following: 

1) If $B'$ contains at least two elements,  then
for all $ b^i, b^j\in B'$, 

$$
b^i\not= b^j\Longleftrightarrow < (b^i_0,b^i_1,\ldots,b^i_h), r>\not=
< (b^j_0,b^j_1,\ldots,b^j_h), r>
$$

2) $f$ doesn't divide $p$ iff $B'\not= \emptyset$, and in this case
there is a unique $k_0$ such that 
$< (b^{k_0}_0,b^{k_0}_1,\ldots,b^{k_0}_h), r>=$ inf$\lbrace < (b^k_0,b^k_1,\ldots,b^h_h), r>; b^k\in B'\rbrace$.}

\end{theorema}

\begin{definicion}  {\rm (see [3])
The integer $< (b^{k_0}_0, b^{k_0}_1,\ldots,b^{k_0}_h), r>$
is called  {\bf formal intersection multiplicity}  
of $p$ with respect to $(r, g)$ and will be 
denoted by  fint$(p, r, g)$.}\end{definicion}

\medskip

\noindent Now we recall the notion of { \bf generalized Newton
polygon}.   Let $p$ be a monic polynomial of ${\bf R}$ of degree $n$
in $y$
and consider a monic polynomial $q$ of ${\bf R}$ of degree  
${\displaystyle\frac{n}{d}}$ in $y$, where
$d$ is a divisor  of  $n$. Let

$$p=q^{d}+\alpha _1(x,y)q^{d-1}+\ldots+\alpha_{d}(x,y)$$

\noindent be the expansion  of  $p$ with respect to the powers of  $q$, and 
consider the sequences  $r$, $g$ defined above.
One associates with  $p$ the generalized Newton polygon  which is 
defined as the union of all compact sides of the convex hull in ${\bf R^2}$  of the set formed by the points
(fint$(\alpha_k, r, g), (d-k).$fint$(q, r, g))$  for all $1\leq k\leq d$. It will be denoted by
GNP$(p,q,r,g)$ (see [3]). \newline
With these notations one has the following:

\medskip

\noindent {\bf Irreducibility criterion} (see [3])

\medskip

\noindent Write  $p= y^n+a_1(x)y^{n-1}+\ldots+a_n(x) \in {\bf R}$ and
assume, possibly after a change of variables, that $a_1(x)=0$. Consider 
the sequences  $r_k,g_k,d_k$ defined in the following way:

 $r_0=d_1=n$

$g_1=y, r_1=$ int$(p,g_1)$, and for all  $k\geq 2$:

$d_{k}=$ gcd $(r_0,r_1,\ldots,r_{k-1}), g_{k}=$ App$_{d_{k}}(p), r_{k}=$
int$(p,g_k)$.

\noindent $p$ is irreducible if and only if the following conditions hold.

1) There is $h\in N$ such that  $d_{h+1}=1$.

2) For all $k=1,\ldots,h-1$, $r_{k+1}d_{k+1} > r_kd_k$.

3) Set $p=g_{h+1}$ and let for all $k=1,\ldots,h$, 
$e_k={\displaystyle\frac{d_k}{d_{k+1}}}$ and 
$g^k=(g_1,\ldots,g_k)$. Let also for all $k=1,\ldots,h$,
${\displaystyle\frac{r}{d_{k+1}}}=({\displaystyle\frac{r_0}{d_{k+1}}},
\ldots,{\displaystyle\frac{r_k}{d_{k+1}}})$. Then for all $k=1,\ldots, h$, the
generalized Newton polygon GNP$(g_{k+1},g_k,
{\displaystyle\frac{r}{d_{k+1}}}, g^k)$ is the line
segment joining $(0,{\displaystyle\frac{r_k}{d_{k+1}}}.e_k)$ and
$({\displaystyle\frac{r_k}{d_{k+1}}}.e_k,0)$.

\medskip

\noindent {\bf Remarks 2.3:} i) (see [3]) Suppose that $p$ is irreducible, and 
let  $r=(r_0,r_1,\ldots,r_h)$ and $g =(g_1=y,
g_2,\ldots,g_h,g_{h+1}=p)$ be the sequences defined above. Let 
$p'$ by a polynomial of ${\bf R}$ and consider the expansion of
$p'$ with respect to the sequences $r,g$ (see Theorem 2.1.). If the 
corresponding set  $B'$ is non empty, then  fint$(p',r,g)=$ int$(p',p)$.

\vskip0.1cm

\noindent ii) Part 3) of the criterion can be precised as follows:  
whenever  $p$ is irreducible, 
the generalized Newton polygon GNP$(g_{k+1},g_k,{\displaystyle\frac{r}{d_{k+1}}},g^k)$ 
just contains the two extremal points for all $k=1,\ldots,h-1$. In 
fact let $k\in\lbrace 1,\ldots,h\rbrace$ and 
let $g_{k+1}=g_k^{e_k}+\alpha_2(x,y).g_k^{e_k-2}+\ldots+\alpha_{e_k}(x,y)$ 
be the expansion of $g_{k+1}$ w.r.t. $g_k$, then we have:

a) int$(g_{k+1},\alpha_{e_k}(x,y))={\displaystyle\frac{r_k}{d_{k+1}}}.e_k=$ int$(g_{k+1},g_k^{e_k})$.

b) For all $i=2,\ldots,e_{k}-1$, int$(g_{k+1},\alpha_{i}(x,y)) > {\displaystyle\frac{r_k}{d_{k+1}}}.i$,  
In particular int$(g_{k+1},\alpha_{i}(x,y))+{\rm int}(g_{k+1},g_k^{e_k-i})>{\displaystyle\frac{r_k}{d_{k+1}}}.e_k$.

\vskip0.1cm

\noindent iii) (see [13]) As an immediate consequence of ii) a) we have that, 
for all $k=1,\dots, h$,
$$ {\displaystyle\frac{r_k}{d_{k+1}}}.e_k\in 
e_k.\Gamma(g_{k})=<{\displaystyle\frac{r_0}{d_{k+1}}},\dots, 
{\displaystyle\frac{r_{k-1}}{d_{k+1}}}>.$$
In particular 
$$ r_k.e_k\in <r_0,\dots, r_{k-1}>.$$

\medskip

\section{Constructing the equisingularity class}

\noindent In this Section we fix a semigroup of nonnegative integers
$\Gamma=<r_0,r_1,\ldots,r_h>$, and  we set $d_1=r_0$ and
$d_{k+1}=$ gcd$(r_0,\ldots,r_k)$ for all $k=1,\ldots,h$ (we set by
convention $r_{h+1}=d_{h+2}=+\infty$). Moreover  we assume that  
$d_{k+1}=1$ and that $r_{k+1}.d_{k+1}> r_k.d_k$ for all
$k=1,\ldots,h$ (*) (this condition appeared in the irreducibility
criterion in Section 2). This implies that the
sequence $(r_1,\ldots,r_h)$ is strictly increasing. This also holds if
$r_1$ is replaced by $r_0$.

Let $\tilde{\bf R}$ denote the set of all irreducible monic polynomials $f$ 
of ${\bf R}$ of the form $f=f(x,y)=y^n+a_2(x).y^{n-2}+\ldots+a_n(x)$. Condition $(*)$ implies that there exists a
polynomial $f\in \tilde{\bf R}$ such that
$\Gamma=\Gamma(f)$ (see [13]). Here we give the generic forms of all these
polynomials, i.e., we describe the set of elements of $\tilde{\bf R}$
having the semigroup $\Gamma$.  The construction can be performed with 
respect to the arrangements $(r_0,r_1,r_2,\ldots,r_h)$ and
$(r_1,r_0,r_2,\ldots,r_h)$. We shall perform it with respect to the
first arrangement. The polynomials that we would get with respect to 
the second arrangement are those obtained by  exchanging $x$ and $y$.

\medskip
\noindent In this and in the following Sections  we shall assume that 
$r_0,\dots, r_h$ form a minimal system of generators for $\Gamma$.
This condition can be reformulated equivalently as a numerical criterion. This
 is what we are going to do next. First we derive a 
useful identity: Set for all $1\leq i\leq h, e_i={\displaystyle\frac{d_i}{d_{i+1}}}$. 
For all $2\leq k \leq h$ we have:

\begin{eqnarray*}\displaystyle\sum_{i=1}^{k-1} (e_i -1)r_i&=&
\displaystyle\sum_{i=1}^{k-1} (r_{i+1}- (m_{i+1} -m_i)) -\displaystyle\sum_{i=1}^{k-1} r_i\\
&=&r_k-r_1-\displaystyle\sum_{i=1}^{k-1}(m_{i+1}-m_i)
= r_k-m_k\quad(**)
\end{eqnarray*}

\noindent Now we are ready to prove the following 

\begin{lema}  Suppose that $r_0,\dots, r_h$ satisfy condition $(*)$.  These 
numbers form a minimal system of generators for the 
semigroup $\Gamma$ if and only if $d_1>d_2>\cdots>d_h>d_{h+1}$.\end{lema}

\begin{demostracion}{.}  Remark that in general $d_k\geq d_{k+1}$ for 
all $k=1,\dots, h$. Moreover, recall that the minimality of the system 
of generators is equivalent to the condition that 
$r_{k}\notin <r_0,\dots, r_{k-1}>$ for all $k=1,\dots, h$. First 
suppose that
 this condition is not fulfilled  for some index $k$. Then  
 $d_{k}=$gcd$(r_0,\dots,r_{k-1})$ certainly divides $r_{k}$. 
Hence $d_{k+1}=$gcd$(d_{k},r_{k})=d_{k}$. For the converse fix an 
index $k$ and suppose  that $d_{k+1}=d_k$.  Then there exist some 
integers $\alpha_0,\dots,\alpha_{k-1}$ such that 

$$
r_k = \displaystyle\sum_{i=0}^{k-1} \alpha_ir_{i}.\qquad\qquad(1)
$$ 

\noindent Now let $\beta_{k-1}$ be the (nonnegative) remainder of the euclidean
division of $\alpha_{k-1}$ by 
$e_{k-1}={\displaystyle{{d_{k-1}}\over{d_{k}}}}$. Since 
the semigroup
$\Gamma$ verifies condition (*),  $\Gamma$ is the semigroup of a
polynomial of $\tilde{\bf R}$, in particular, by Remarks 2.3 iii), 
$r_{k-1}.(e_{k-1})\in<r_0,\dots, r_{k-2}>$; hence we can transform  $(1)$ in such a way that $\alpha_{k-1}=\beta_{k-1}$. 
If we successively perform the same procedure  for the indices  $k-2,\dots,1$, we finally obtain that in  (1) 
$0\leq\alpha_i< e_i$ for all $i=1,\dots,k-1$. Now by $(**)$  

$$
\alpha_0r_0 = r_k-\displaystyle\sum_{i=1}^{k-1} \alpha_ir_{i}>r_k-\displaystyle\sum_{i=1}^{k-1} (e_i-1)r_i=m_k>0,
$$

\noindent hence $\alpha_0>0$. This proves that $r_{k}\in <r_0,\dots, r_{k-1}>$ and completes the proof.\end{demostracion}

\noindent The construction of the generic form of all polynomials $f\in\tilde{\bf R}$ having $\Gamma$ as a semigroup is based on the notion of {\bf generalized
Newton polygons} introduced in Section 2. We shall recursively
construct the sequence of approximate roots $g_1,\ldots,g_h,g_{h+1}=f$.

\medskip


\noindent Let $g_1=$ App$_{d_1}(f)$ (and recall that, since
$a_1(x)=0$, then $g_1=y$). From Section 2 we know that $g_2=$
App$_{d_2}(f)$ satisfies:

i) $\Gamma(g_2)=<{\displaystyle\frac{r_0}{d_2}},{\displaystyle\frac{r_1}{d_2}}>$.

ii)  GNP$(g_2,g_1,{\displaystyle\frac{r}{d_2}}, g^1)$ is the line segment 
joining the two points $(0,{\displaystyle\frac{r_1}{d_2}}.e_1)$ 
and $({\displaystyle\frac{r_1}{d_2}}.e_1,0)$. 

\noindent By virtue of  part ii) of Remarks 2.3, this  yields  to the following generic form of $g_2$:

$$
g_2 = y^{{\frac{r_0}{d_2}}} +
a.x^{{\frac{r_1}{d_2}}}+
\sum_{i.{\frac{r_0}{d_2}}+
j.{\frac{r_1}{d_2}} >
{\frac{n}{d_2}}.{\frac{r_1}{d_2}}; 0\leq j<
{\frac{d_1}{d_2}}={\frac{r_0}{d_2}}} a_{ij}x^iy^j,
$$

\noindent where $a\in {\bf K}-0$ and for all
$(i,j)$, $a_{ij}\in {\bf K}$. 

\noindent Suppose that we have the generic forms of $g_1,\ldots,g_k$
and consider the expansion of $g_{k+1}$ with respect to $g_k$:

$g_{k+1} =
g_k^{e_k}+\alpha_2(x,y)g_k^{e_k-2}+\ldots+\alpha_{e_k}(x,y)$.

\noindent From Section 2 we know that:

i) $\Gamma(g_{k+1})=<{\displaystyle\frac{r_0}{d_{k+1}}},\ldots,
{\displaystyle\frac{r_k}{d_{k+1}}}>$.

ii) GNP$(g_{k+1},g_k,{\displaystyle\frac{r}{d_{k+1}}},g^k)$ is the 
line segment joining the two points $(0,{\displaystyle\frac{r_k}{d_{k+1}}}.
e_k)$ and $({\displaystyle\frac{r_k}{d_{k+1}}}.e_k,0)$.

\noindent It follows from Remarks 2.3 that 

(1) int$(g_{k+1},\alpha_{e_k}(x,y))={\rm
fint}(\alpha_{e_k}(x,y),{\displaystyle\frac{r}{d_{k+1}}}, g^k)=
{\displaystyle\frac{r_k}{d_{k+1}}}.e_k$,

\noindent and that for all $i=2,\ldots,e_k-1$:

(2) int$(g_{k+1},\alpha_i(x,y))={\rm
fint}(\alpha_i(x,y),{\displaystyle\frac{r}{d_{k+1}}},g^k)>
{\displaystyle\frac{r_k}{d_{k+1}}}.i$,

\noindent Let $B^k$ be the set of all 
$\theta=(\theta_0,\ldots,\theta_{k-1})\in {\bf N}^k$ such that 
 $ 0\leq \theta_j < e_j$ for all $j=1\ldots,k-1$, then
associate with all $\theta\in B^k$ the ``monomial''
$M_{\theta}=x^{\theta_0}.g_1^{\theta_1}.\ldots.g_{k-1}^{\theta_{k-1}}$.
For all $i\in {\bf N}$ and for all $\theta\in B^k$, we say that $M_{\theta}$ is of type $(k,i,1)$ 
(resp. of type $(k,i,2)$) if  

$$
{\displaystyle\frac{r_k}{d_{k+1}}}.i=\theta_0.{\displaystyle\frac{r_0}{d_{k+1}}}+
\theta_1.{\displaystyle\frac{r_1}{d_{k+1}}}
+\ldots +\theta_{k-1}.{\displaystyle\frac{r_{k-1}}{d_{k+1}}}
$$ 

\noindent (resp. 

$$
{\displaystyle\frac{r_k}{d_{k+1}}}.i<\theta_0.{\displaystyle\frac{r_0}{d_{k+1}}}
+\theta_1.{\displaystyle\frac{r_1}{d_{k+1}}}
+\ldots +\theta_{k-1}.{\displaystyle\frac{r_{k-1}}{d_{k+1}}}).
$$ 

\noindent Let 
$E(k,i,1)$ (resp. $E(k,i,2)$)
be the set of monomials $M_{\theta},\theta\in B^k$, of type $(k,i,1)$ (resp. of type $(k,i,2)$). Since 
${\displaystyle\frac{r_k}{d_{k+1}}}.e_k\in <{\displaystyle\frac{r_0}{d_{k+1}}},\ldots,
{\displaystyle\frac{r_{k-1}}{d_{k+1}}}>$,
then  $E(k,e_k,1)$  is reduced to one element. If we write this element as
$M_{\theta^k}=x^{\theta_0}.g_1^{\theta_1}.\ldots.g_{k-1}^{\theta_{k-1}}$, then 
$(\theta_0,\theta_1,\ldots,\theta_{k-1})$ can be
calculated by euclidean division.

\noindent Using Remark 2.3., (1) and (2) this leads to the following generic forms
of $\alpha_2,\ldots,\alpha_{e_k}:$

$$
\alpha_{e_k}=a.M_{\theta^k}+\sum_{M_{\theta}\in
E(k,e_k,2)}a_{\theta}.M_{\theta},
$$

\noindent (resp. for all $i=2,\ldots, e_k-1$,

$$
\alpha_{i}=\sum_{M^i_{\theta}\in E(k,i,2)}a^i_{\theta}.M^i_{\theta}),
$$

\noindent where $a\in {\bf K}-0$, and for all $\theta,
a_{\theta}\in {\bf K}$ (resp. for all $\theta$ and for
all $i=2,\ldots, e_k-1$, $a^i_{\theta}\in {\bf K}$).

\medskip

\begin{nota}{\rm  We proved that, if $\Gamma$ is the semigroup of a polynomial 
$f\in\tilde{\bf R}$, then $f$ and its approximate roots $g_1,\ldots,g_h$ 
belong to the set of polynomials constructed above. Conversely, let $(g_1,\ldots,
g_h,g_{h+1}=f)$ be as above, then part ``only if'' of the irreducibility 
criterion of Abhyankar shows that $f$ is irreducible.}\end{nota}

\begin{nota} {\rm Given $1\leq k \leq h$, it follows from the above
construction that a polynomial $g_{k+1}$ may have an infinite number 
of monomials. In particular the  above construction is not algorithmic.
Remark however that $g_{k+1}$ is obtained from the sum 
$g_{k}^{e_k}+a.M_{\theta^{k}}, a\in {\bf K}-0$ by adding monomials 
that verify some conditions. This suggests the introduction of the
following set of polynomials: let $G_1=y$ and for all $1\leq k\leq h$,
$G_{k+1}=G_k^{e_k}-M_{{\theta}^k}$. For all $1\leq k\leq h$, $g_k$ is obtained
from $G_k$ in an obvious way. Set $G=(G_1,\ldots,G_h,G_{h+1})$ and call it
the canonical element of the set of all $(g_1,\ldots,g_h,g_{h+1})$
constructed above. The above calculation leads to an algorithm that
computes this canonical element. It is based on euclidean division
in ${\bf N}$. The different steps can be summarized as follows:

\medskip

i) Consider a sequence of integers $r_0 <\ldots <r_h$.

ii) Compute the gcd sequence $d=(d_1,\ldots,d_h,d_{h+1})$ such
that $d_1=r_0$ and for all $2\leq k \leq h+1, d_{k}={\rm gcd}(r_{k-1},d_{k-1})$. 
 Let $e_k={\displaystyle\frac{d_k}{d_{k+1}}}$ for all $1\leq k\leq h$. 

iii) If either $d_{h+1} > 1$, or $r_k\cdot e_k\geq r_{k+1}$ for
at least one $1\leq k \leq h$, then the sequence $(r_0, \ldots, r_h)$ is
not the semigroup of an irreducible polynomial of $\tilde{\bf R}$.

iv) Assume that  $d_{h+1}=1$ and that $r_k\cdot d_k< r_{k+1}\cdot d_{k+1}$ for
all $1\leq k \leq h-1$.
   
   \quad\quad a) If $d_k=d_{k+1}$ for some $1\leq k\leq h$, then eliminate 
   $r_{k}$ from the $r$-sequence of i).
   
   \quad\quad b) Assume that $d_1 > d_2 >\ldots >d_h > d_{h+1}=1$. Then 
   for all $1\leq k \leq h$, compute (the unique) $\theta^k=(\theta_0^k,\ldots,
   \theta_{k-1}^k)$ such that $ 0\leq \theta_j^k < e_j$ 
   for all $j=1\ldots,k-1$ and 
   ${\displaystyle\frac{r_k}{d_{k+1}}}.e_k=\theta_0^k.
   {\displaystyle\frac{r_0}{d_{k+1}}}+
\theta_1^k.{\displaystyle\frac{r_1}{d_{k+1}}}
   +\ldots +\theta_{k-1}^k\cdot {\displaystyle\frac{r_{k-1}}{d_{k+1}}}$.
   
   \quad\quad c) The canonical element is $G=(G_1,\ldots,G_h,G_{h+1})$
   where $G_1=y$ and for all $2\leq k\leq h+1$, 
   $G_k=G_{k-1}^{e_{k-1}}-x^{\theta_0^k}y^{\theta_1^k}
   \cdot\ldots\cdot G_{k-1}^{\theta_{k-1}^k}$.}
 \end{nota}

\noindent This algorithm has been implemented with {\it Mathematica} 
(see [8]), and {\it Maple}: the input is an increasing sequence of positive integers.
  Then the output is "false" if  this sequence 
 does not generate the semigroup of an irreducible polynomial of $\tilde{\bf R}$. Otherwise,
we get  the canonical element described above.

\medskip

\noindent Note that our implementation is based on the following: 
given $r_0,r_1,\dots,r_{k-1}$, we need to  compute the unique 
 $\theta^k=(\theta_0^k,\ldots,
   \theta_{k-1}^k)$ such that $ 0\leq \theta_j^k < e_j$ 
   for all $j=1\ldots,k-1$ and 
 ${\displaystyle\frac{r_k}{d_{k+1}}}.e_k=\theta_0^k.
   {\displaystyle\frac{r_0}{d_{k+1}}}+
\theta_1^k.{\displaystyle{\frac{r_1}{d_{k+1}}}}+\cdots +\theta_{k-1}^k.
{\displaystyle{\frac{r_{k-1}}{d_{k+1}}}}$. Instead 
of applying the Euclidean division, 
We have preferred to  scann lists of values, namely the set
of values $(a_0,a_1,\ldots,a_{k-1})$ where for all $i\geq 1,
0\leq a_i < e_i$ and $0\leq a_0 \leq
{\displaystyle{\frac{r_k}{d_{k+1}}}.e_k.{\frac{d_{k+1}}{r_0}}=
{\frac{r_k.e_k}{r_0}}}$. The cardinality of this set is:

$$
{\frac{r_k.e_k}{r_0}}.\prod_{i=1}^{k-1}e_k={\frac{r_k.e_k}{r_0}}.{\frac{d_1}{d_k}}={\frac{r_k}{d_{k+1}}}
$$

\noindent In conclusion the set of the values scanned in the algorithm
is
bounded by ${\displaystyle{\sum_{k=1}^h {\frac{r_k}{d_{k+1}}}}}$.

\begin{nota} {\rm An element $f$ whose semigroup is $\Gamma$ can also be
calculated by using the theory of Gr\"obner bases: a reduced Gr\"obner basis
with respect to  any well-ordering on $\bf N^3$ that 
eliminates $t$ from the equations $x-t^n, y-t^{m_1}-\ldots-t^{m_r}$ contains
a unique polynomial $f(x,y)$. If we consider $f$ as an element of ${\bf K}[[x,y]]$, 
then obviously $\Gamma=<r_0,\ldots,r_h>$ is the semigroup of $f$. It is
well known that the complexity of a Gr\"obner basis is in general
doubly exponential. Moreover, the algorithm computes more than we
need. We
think that our option is more natural in view of 
our situation, especially because of its complexity
and that the output is expressed in terms of the polynomial $f$.}
\end{nota}

\noindent {\bf Example:}

\noindent Let  $\Gamma = <8,12,50,101>$. Here $h=3$, the $r$-sequence is
$r=(8,12,50,101)$, and the gcd-sequence is $d=(8,4,2,1)$. Moreover, 
$e_1=e_2=e_3=2$. Let us construct
the canonical element $G=(G_1,G_2,G_3,G_4)$ following the algorithm 
above. Here we start directly from  point iv), b): 

$k=1: {\displaystyle\frac{r_1}{d_2}}\cdot{e_1}=3\cdot 2= 
\theta_0^1\cdot{\displaystyle\frac{r_0}{d_2}}=\theta_0^1\cdot 2$ implies
that $\theta_0^1=3$.

$k=2:  50={\displaystyle\frac{r_2}{d_3}}\cdot {e_2}=\theta_0^2\cdot
{\displaystyle\frac{r_0}{d_{3}}}+
\theta_1^2\cdot {\displaystyle\frac{r_1}{d_{3}}}= \theta_0^2\cdot 4+
\theta_1^2\cdot 6$ with $0\leq \theta_1^2 <2$.
This implies that $\theta_1^2=1$, and $\theta_0^2=11$.

$k=3: 202={\displaystyle\frac{r_3}{d_4}}\cdot {e_3}=\theta_0^3\cdot
{\displaystyle\frac{r_0}{d_{4}}}+
\theta_1^3\cdot {\displaystyle\frac{r_1}{d_{4}}}+\theta_2^3\cdot
{\displaystyle\frac{r_2}{d_{4}}}=\theta_0^3\cdot 8+
\theta_1^3\cdot 12+\theta_2^3\cdot 50$ with $0\leq \theta_1^3, \theta_2^3<2$.
This implies that $\theta_2^3=1, \theta_1^3=0, \theta_0^3=19$.

\medskip

\noindent In particular, $G_1=y, G_2=G_1^2-x^3=y^2-x^3, 
G_3=G_2^2-x^{11}\cdot G_1=(y^2-x^3)^2-x^{11}\cdot y, 
G_4=G_3^2-x^{19}\cdot G_2= [(y^2-x^3)^2-x^{11}\cdot y]^2-
x^{19}\cdot (y^2-x^3).$

\medskip

\noindent  With the same notations as above, the
set  of elements $(g_1, g_2, g_3, g_4=f)$ is then given by:

$$g_1=y.$$ 

$$
g_2= y^2+\alpha_2(x)= y^2+ax^3+\sum_{M_{\theta}\in
E(1,2,2)}a_{\theta}.M_{\theta},
$$ 

\noindent where $a\in {\bf K}-0$, and for all $\theta$, one has
 $a_{\theta}\in {\bf K}$ and $M_{\theta}=x^{\theta_0}$, with $6<2\theta_0$.
Moreover 

$$g_3= g_2^2 +\alpha'_2(x,y)= g_2^2 + a' x^{11}y+\sum_{M'_{\theta}\in
E(2,2,2)}a'_{\theta}.M'_{\theta},$$ where

-  $a'\in {\bf K}-0$, and for all $\theta$,  $a'_{\theta}\in {\bf K}$;

-  for all $\theta$, $M'_{\theta}=x^{\theta'_0}y^{\theta'_1}$, with $50 <4\theta'_0 +6\theta'_1$. 

\noindent Finally, 
$$f=g_3^2+\alpha''_2(x,y)= g_3^2 + a'' x^{19}g_2+\sum_{M''_{\theta}\in
E(3,2,2)}a''_{\theta}.M''_{\theta},$$ where

-  $a''\in {\bf K}-0$, and for all $\theta$,  $a''_{\theta}\in {\bf K}$;

-   for all $\theta$, $M''_{\theta}=x^{\theta_0}y^{\theta_1}g_2^{\theta_2}$, with $202<8\theta''_0 +12\theta''_1+50\theta''_2$. 
Hence the generic form of all polynomials having $\Gamma$ as a semigroup is the following:

$$f= [ (y^2+ax^3+F)^2+a'x^{11}y+F']^2 +
a''x^{19}(y^2+ax^3+F) +F'',$$
where $a,a',a''\in {\bf K}-0$ and $F$, $F'$ and $F''$ are arbitrary linear combinations of 
monomials from $E(1,2,2)$, $E(2,2,2)$ and  $E(3,2,2)$ respectively.

\begin{nota}{\rm i) The construction above does not depend on the choice 
of the coefficients in the field $\KK$ -provided that it is of 
characteristic zero-, in particular the algorithm described allows
us to work over anysubring $A$ of $\KK$. If 
$A=k[t_1,\ldots,t_m]$ is a polynomial ring over a field $k$ of
characteristic zero and $\KK$ is the algebraic closure of $A$
into its fractions field, then we get the equisingularity
class of the $(t_1,\ldots,t_m)$-generic section.

\medskip

ii) The restriction to the zero characteristic is made
only because of the use of the approximate roots in
the algorithm. If the characteristic of $\KK$ does not
divide $r_0$, then everything above applies -see Remark 1.6.-.
Note that a more general irreducibility criterion has been given
by A. Granja (see [12]), but it does not seem to be in computational
form.}
\end{nota}

\section{ Equisingularity classes with a given Milnor number}

\medskip

\noindent In this Section we generalize the results of
Section 3 in the following way: let $m\in {\bf N}$ be a fixed
integer. If $m\in 2.{\bf N}$, then there exists a
polynomial $f=y^n+a_2(x).y^{n-2}+\ldots+a_n(x)\in\tilde{\bf R}$
such that int$(f_x,f_y)=m$. Here we shall give the generic forms of all
these polynomials. Remark that if $g$ is another 
polynomial of $\tilde{\bf R}$, then $\Gamma(f)=\Gamma(g)$ implies that
int$(f_x,f_y)={\rm int}(g_x,g_y)$. Thus the set of
$f=y^n+a_2(x).y^{n-2}+\ldots+a_n(x)\in \tilde{\bf R}$ such that
int$(f_x,f_y)=m$ is the union  of equisingularity classes. We shall
first prove that this union is finite. This is an
immediate application of the next Proposition. We recall that if a subsemigroup of $\ZZ$ is minimally generated by $h+1$ elements, then $h$ is called the {\it length} of the semigroup.  

\begin{proposicion} \label{majoration} Let $h\in \NN$ and consider a
polynomial $f\in \tilde{{\bf R}}$ such that $h$ is the length of
$\Gamma(f)$. Let $\mu_{h+1}={\rm int}(f_x,f_y)$, and let $r_h$ be
the last generator of $\Gamma(f)$. We have the following:

i) $h=1$ implies that $r_h \geq 3$ and $\mu_{h+1}\geq 2$.  
   
ii) $h=2$ implies that $r_h\geq 13$ and $m_{h+1}\geq 16$.

iii)  More  generally we have:

\noindent 1) $r_h \geq 12.4^{h-2}+\sum_{i=0}^{h-2}4^i={\displaystyle\frac{5}{3}}.2^{2h-1}-{\displaystyle\frac{1}{3}}$.

\noindent 2) $\mu_{h+1}\geq 2+2.\sum_{i=0}^{h-2}4^i+12.\sum_{i=h-2}^{2h-4}2^i={\displaystyle\frac{5}{3}}.2^{2h}-3.2^h+
{\displaystyle\frac{4}{3}}$, assuming that the summation over
negative exponents is $0$.\end{proposicion} 
                
\begin{demostracion}{.} i) In this case,
by Lemma 1.8., $\mu_1=(r_0-1).(r_1-1)$. Furthermore $r_1\geq 2$  and $r_0\geq 2$; otherwise
$\Gamma(f)=<1>$, and $h=0$. On the other hand, gcd$(r_0,r_1)=1$. This
proves that max$(r_0,r_1)=r_1\geq 3$ and $\mu_2\geq 2$. Then our assertion follows. Remark that $r_1=3$ and 
$\mu_2=2$ holds for $f=y^2+ax^3$, where $a\in \KK-0$.

ii) Let $g_2$ be the second approximate root of $f$, then $\mu_3=
d_2.${\rm int}$(g_{2_x},g_{2_y})+(d_2-1)(r_2-1)$. It follows from i)
that ${\displaystyle\frac{r_1}{d_2}}\geq 3$ and that {\rm
int}$(g_{2_x},g_{2_y})=({\displaystyle\frac{r_0}{d_2}}-1)({\displaystyle\frac{r_1}{d_2}}-1)\geq 2$, 
and also that ${\displaystyle\frac{r_0}{d_2}}+{\displaystyle\frac{r_1}{d_2}}
\geq 5$. In particular 

$$
(r_0-d_2)({\displaystyle\frac{r_1}{d_2}}-1)\geq 2.d_2
$$ 

\noindent and 

$$
r_0+r_1\geq 5.d_2.
$$

\noindent Thus:

$$
r_1.{\displaystyle\frac{r_0}{d_2}}\geq d_2+r_1+r_0\geq 6.d_2\geq 12.
$$

\noindent But $r_2-1 \geq r_1.{\displaystyle\frac{d_1}{d_2}} =r_1.{\displaystyle\frac{r_0}{d_2}}$. Finally
$r_2\geq r_1{\displaystyle\frac{r_0}{d_2}}+1\geq 13$, and $\mu_{h+1}\geq 2.d_2+(r_2-1) \geq  4+12=16$. This 
implies our assertion. Note that the lower 
bounds 13 and  16 is sharp:they are satisfied 
for $f=(y^2+a.x^3)^2+b.x^5y,$ where $a,b\in \KK-0$, 
whose semigroup is $\Gamma=<4,6,13>$. 

iii)  We prove the inequalities by induction on $h$. From i) and ii) both
are satisfied for $1\leq h\leq 2$. Assume that $h\geq 3$ and that the
formulas are true for $h-1$. We first prove the inequality 1): Remark first that $r_h \geq
({\displaystyle\frac{r_{h-1}}{d_h}}).d_{h-1}+1$. The quotient ${\displaystyle\frac{r_{h-1}}{d_h}}$ being the last generator of
$\Gamma(g_h)$ which is of length $h-1$, it follows by induction that
${\displaystyle\frac{r_{h-1}}{d_h}} \geq  12.4^{h-3}+\sum_{i=0}^{h-3}4^i={\displaystyle\frac{5}{3}}.2^{2h-3}-{\displaystyle\frac{1}{3}}$. On the other hand,
$d_{h-1} \geq 4$, thus  $r_h\geq 4.({\displaystyle\frac{5}{3}}.4^{2h-3}-{\displaystyle\frac{1}{3}})+1={\displaystyle\frac{5}{3}}.2^{2h-1}-{\displaystyle\frac{1}{3}}$. This is the required inequality.

\noindent  We now prove the inequality 2): Consider to this end the last approximate root $g_h$ of $f$. We
have: $\mu_{h+1}=d_h.{\rm int}(g_{h_x},g_{h_y})+(d_h-1)(r_h-1)$. But $d_h\geq 2$ and $r_h\geq {\displaystyle\frac{5}{3}}.2^{2h-1}-{\displaystyle\frac{1}{3}}$. On the 
other hand, the length of
$\Gamma(g_h)$ being $h-1$, it follows that 

$$
{\rm int}(g_{h_x},g_{h_y})\geq 2+2.\sum_{i=0}^{h-3}
4^i+12.\sum_{i=h-3}^{2h-6}2^i={\displaystyle\frac{5}{3}}.2^{2h-2}-3.2^{h-1}+{\displaystyle\frac{4}{3}}
$$

\noindent In particular 

$$
\mu_{h+1}\geq {\displaystyle\frac{5}{3}}.2^{2h-2}-3.2^{h-1}+{\displaystyle\frac{4}{3}}+{\displaystyle\frac{5}{3}}.2^{2h-1}-{\displaystyle\frac{1}{3}}-1={\displaystyle\frac{5}{3}}.2^{2h}-3.2^{h}+{\displaystyle\frac{4}{3}}
$$

\noindent  This is the required inequality.
\end{demostracion}

\begin{nota}{\rm The bounds of the above Proposition are sharp. More
precisely, for all $h\geq 1$, there is a polynomial $f_h(x,y)\in \tilde{R}$
such that $h$ is the length of $\Gamma(f)$, and that int$(f_h{_x},f_h{_y})={\displaystyle\frac {5}{3}}.2^{2h}-3.2^h+{\displaystyle\frac {4}{3}}$, and if $r_h$ denotes the last generator of $\Gamma(f)$, then 
$r_h={\displaystyle\frac{5}{3}}.2^{2h-1}-{\displaystyle\frac{1}{3}}$. Consider to this end the semigroup $\Gamma_h$ generated by $r_0=2^h$ and  

$$
r_k=2^{h-k}({\displaystyle\frac{5}{3}}.2^{2k-1}-{\displaystyle\frac{1}{3}})
$$

\noindent for all $1\leq k\leq h$ (equivalently  $r_1=2^{h-1}.3, r_2=2^h.3+2^{h-2},\ldots, 
r_{k+2}=2^{h+k}.3+\sum_{i=1}^{k+1}2^{h+k-2i}$ for all   $1\leq k\leq h-2$). Clearly $r_h={\displaystyle\frac{5}{3}}.2^{2h-1}-{\displaystyle\frac{1}{3}}$, and 
the $d$-sequence is given by $d_k=2^{h+1-k}, 1\leq k\leq h+1$. Furthermore, 
$r_kd_k < r_{k+1}d_{k+1}$  for all 
$1\leq k\leq h$. It follows that $\Gamma_h$ is the semigroup of a 
polynomial of $\tilde{R}$. We shall prove by
induction that the Milnor number of such a polynomial 
is ${\displaystyle\frac {5}{3}}.2^{2h}-3.2^h+{\displaystyle\frac {4}{3}}$. Denote this
 number by $\mu_{h+1}$ and recall that 
 $\mu_{h+1}=\sum_{k=1}^h({\displaystyle\frac {d_k}{d_{k+1}}}-1)r_k-r_0+1$. Since 
 ${\displaystyle\frac {d_k}{d_{k+1}}}=2$ for all $1\leq k\leq h$, then:

$$
\mu_{h+1}=(\sum_{k=1}^hr_k)-r_0+1=\sum_{k=1}^h 2^{h-k}({\displaystyle\frac{5}{3}}.2^{2k-1}-{\displaystyle\frac{1}{3}})-2^h+1
$$

\noindent Which is nothing but ${\displaystyle\frac{5}{3}}.2^{2h}-3.2^h+
{\displaystyle\frac {4}{3}}$. This proves our assertion.}

\end{nota}

\begin{corolario} Let $m\in 2.\NN$, then one can effectively compute  the set of
irreducible polynomials $f\in {\tilde{\bf R}}$ such that $m={\rm int}(f_x,f_y)$.\end{corolario}

\begin{demostracion}{.} It follows from Proposition 4.1. that the
length $h$ of the semigroup of a polynomial 
$f\in {\tilde{\bf R}}$ with $m={\rm int}(f_x,f_y)$ takes a finite
number of values. In fact, easy calculations show that $h$ must verify
the inequality: $2^h\leq M={\displaystyle\frac{9+\sqrt{1+60m}}{10}}$. In particular,
$h\leq{\displaystyle\frac {{\rm ln}(M)}{{\rm ln}(2)}}$. Let $H=\lbrace h\in {\bf N}; h\leq
{\displaystyle\frac{{\rm ln}(M)}{{\rm ln}(2)}}\rbrace$. Given $h\in A$ we shall effectively
construct the set $\Sigma$ of all the sequences
$(r_0,r_1,\ldots,r_h)$ which minimally generate a semigroup of a polynomial
$f\in {\tilde{\bf R}}$ of the required Milnor number. The steps 
of the algorithm can 
be summarized as follows:
\end{demostracion}

\bigskip

\noindent Set $m=\mu_{h+1}$. We want to calculate the set of $(\mu_h,r_h,d_h)$ with
the following equality: 

\vskip0.1cm
(E1) $\mu_{h+1}=\mu_hd_h+(r_h-1)(d_h-1)$

\vskip0.1cm

\noindent Recall that we have the following restrictions:

\vskip0.1cm

i) $d_h\geq 2$

\vskip0.1cm

ii) $r_h \geq {\displaystyle\frac{5}{3}}.2^{2h-1}-{\displaystyle\frac{1}{3}}$.

\vskip0.1cm

iii) gcd$(r_h,d_h)=1$.

\vskip0.1cm

iv) $\mu_h\geq {\displaystyle\frac{5}{3}}.2^{2h-2}-3.2^{h-1}+
{\displaystyle\frac{4}{3}}$ 

\vskip0.1cm

v) $\mu_1=0$, and for all $h\geq 2$, $\mu_h={\displaystyle{\frac{\mu_{h+1}-(d_h-1)(r_h-1)}{d_h}}}$ is an even integer.

\vskip0.1cm

\noindent Now equality (E1) gives $(d_h-1)(r_h-1)=\mu_{h+1}-\mu_hd_h$,
and by iv)

$$
-\mu_hd_h \leq -[{\displaystyle\frac{5}{3}}.2^{2h-2}-3.2^{h-1}+
{\displaystyle\frac{4}{3}}]d_h
$$

\noindent in particular

$$
(r_h-1)(d_h-1)\leq \mu_{h+1}-[{\displaystyle\frac{5}{3}}.2^{2h-2}-3.2^{h-1}+
{\displaystyle\frac{4}{3}}].d_h
$$

\noindent This gives us  the following upper bound for $r_h$:

$$
r_h\leq {\mu_{h+1}\over{d_h-1}}-[{\displaystyle\frac{5}{3}}.2^{2h-2}-3.2^{h-1}+
{\displaystyle\frac{4}{3}}].{d_h\over {d_h-1}}+1
$$


\begin{corolario} The above equality  with (ii) give:

$$
(E2)\quad {\displaystyle\frac{5}{3}}.2^{2h-1}-{\displaystyle\frac{1}{3}}\leq 
r_h\leq {\mu_{h+1}\over{d_h-1}}-[{\displaystyle\frac{5}{3}}.2^{2h-2}-3.2^{h-1}+
{\displaystyle\frac{4}{3}}].{d_h\over {d_h-1}}+1
$$

\end{corolario}

\medskip

\noindent In the following we shall refine the lower bound 
of Corollary 4.4. We start with the following technical Lemma:


\begin{lema} $\sum_{i=1}^h(e_i-1)r_i=r_hd_h-m_h$. In particular,
$\mu_{h+1}=\sum_{i=1}^h(e_i-1)r_i-r_0+1=r_hd_h-m_h-r_0+1$, where we recall that 
$e_i={\displaystyle\frac{d_i}{d_{i+1}}}$ for all $1\leq i\leq h$.

\end{lema}

\begin{demostracion}{.} Applying identity (**) of Section 3 with $k=h$ we get:

$$
\sum_{i=1}^{h-1}(e_i-1)r_i=r_h-m_h
$$

\noindent Now adding $(e_h-1)r_h=(d_h-1)r_h$ to the equality we get our assertion.
\end{demostracion}

\noindent Lemma 4.5. with equality (E1) imply that $r_hd_h=\mu_{h+1}+m_h+r_0-1$. On the other hand,
we have, with the notation $m_0=r_0$, that for all $1\leq k \leq h$, $m_k-m_{k-1}\geq d_{k+1}$.
Adding these inequalities we get 

$$
\mu_h\geq m_0+d_2+\ldots+d_h+d_{h+1}=d_1+d_2+\ldots+d_h+1.
$$

\noindent But for all $1\leq k\leq h, d_k\geq 2^{h-k}.d_h$, so $\mu_h\geq d_h.(2^h-1)+1$. Since
$r_0=d_1\geq 2^{h-1}.d_h$, we get the following:

$$
r_h \geq \quad{\rm max}({\frac{5}{3}}.2^{2h-1}-{\frac{1}{3}},{\mu_{h+1}\over d_h}+(3.2^{h-1}-1))
$$


\noindent Now equality (E1) implies that 
${\displaystyle {\frac{\mu_{h+1}}{d_h}}}=\mu_h+(r_h-1)(1-{\displaystyle{\frac{1}{d_h}}})$. But
$d_h\geq 2$, thus, using inequalities of Proposition 4.1. we get:

 $$
{\frac{\mu_{h+1}}{d_h}}\geq ({\frac{5}{3}}.2^{2h-2}-3.2^{h-1}+{\frac{4}{3}})+
({\frac{5}{3}}.2^{2h-1}-{\frac{4}{3}}).{\frac{1}{2}}={\frac{5}{3}}.2^{2h-1}-3.2^{h-1}+{\frac{2}{3}}
$$

\noindent In particular ${\displaystyle {\rm max}({\frac{5}{3}}.2^{2h-1}-{\frac{1}{3}},{\mu_{h+1}\over d_h}+(3.2^{h-1}-1))={\mu_{h+1}\over d_h}+(3.2^{h-1}-1)}$. This implies the following:

$$
(E3)\quad {\frac{\mu_{h+1}}{d_h}}+(3.2^{h-1}-1)\leq r_h \leq 
{\frac{\mu_{h+1}}{d_h-1}}-[{\frac{5}{3}}.2^{2h-2}-3.2^{h-1}+{\frac{4}{3}}].{\frac{d_h}{d_{h}-1}}+1
$$

\noindent We shall now use inequality (E3) in order to give an upper bound for $d_h$ (a lower bound
being $2$). Remark to this end that ${\displaystyle {\frac{\mu_{h+1}}{d_h-1}}-[{\frac{5}{3}}.2^{2h-2}-3.2^{h-1}+{\frac{4}{3}}].{\frac{d_h}{d_{h}-1}}+1-({\frac{\mu_{h+1}}{d_h}}+(3.2^{h-1}-1))\geq 0}$.
If we set 
$p=({\displaystyle{\frac{5}{3}}.2^{2h-2}-3.2^{h-1}+{\frac{4}{3}}})$ and $q=3.2^{h-1}-2$, then
an obvious analysis of the above inequality shows that it is equivalent to say that $(p+q).d_h^2-qd_h-\mu_{h+1}\leq 0$, which is true if and only if the following holds:

$$
\mbox {(E4)  } 2\leq d_h\leq {\frac {q+\sqrt{q^2+4\mu_{h+1}.(p+q)}}{2.(p+q)}}=
{\frac {3.2^{h-1}-2+\sqrt{(3.2^{h-1}-2)^2+4\mu_{h+1}.({\frac{5}{3}}.2^{2h-2}-{\frac{2}{3}})}}{{\frac{10}{3}}.2^{2h-2}-{\frac{4}{3}}}}
$$

\noindent {\bf The algorithm}: The two integers $\mu_{h+1}$ and $h$ being fixed, inequality (E4) 
determines the set $D^h$ of possible
values of $d_h$. Each value of $d_h$ gives rise, using inequality (E3), to a set -denoted 
$R^h_{d_h}$- of possible
values of $r_h$ (Remark that ${\displaystyle{\frac{\mu_{h+1}-(d_h-1)(r_h-1)}{d_h}}}$ 
should be an even integer). We get this way the set -denoted $P^h_{d_h}$- of possible values of $(\mu_h,r_h,d_h)$. Now we restart with the set of $\mu_h$... This procees shall stop constructing a
set of lists of length $h$.  The set of semigroups corresponding to $\mu_{h+1}$ is a subset of this list and can be easily calculated. Remark that if $h=1$, then $\mu_1=0$ 
and $\mu_2=(r_1-1)(d_1-1)$ by condition v). In 
this case, the values  of $(r_1,d_1=r_0)$ can 
also be obtained from the set of divisors of 
$\mu_2$.

\begin{exemple}{\rm We perform an explicit computation for $\mu_{h+1}=28$. In this
case, \break $M={\displaystyle\frac{9+\sqrt{1+60.28}}{10}}=5$, so $H=\lbrace
h; 1\leq h \leq {\displaystyle\frac{{\rm ln}(5)}{{\rm ln}(2)}}\rbrace=
\lbrace 1,2\rbrace$. 

\par\smallskip

\noindent 1) $h=1$: In this case, since $28=1*28=2*14=4*7$, 
then $(r_1,d_1)\in \lbrace (2,29),(3,15),(5,8)\rbrace$ and 
condition iii) eliminates  $(3,15)$. We get this way the 
semigroups  $<2,29>$ and $<5,8>$. The  
canonical representative of the equisingularity class of the first
one (resp. the second one) is  $y^2-x^{29}$ (resp. $y^5-x^{8}$).

\par\smallskip

\noindent 2) $h=2$: Inequality (E4) implies in this case 
that $2\leq d_2\leq {\displaystyle{\frac{4+\sqrt{688}}{12}}<3}$. In particular 
$D^2=\lbrace 2\rbrace$. 

\par\smallskip

\noindent Now inequality (E3) implies that 
${\displaystyle{\frac{28}{2}}}+5=19\leq r_2\leq 28-4+1=25$, and 
with conditions iii), v), we get $R^2_{2}=\lbrace 21,25\rbrace$. If $r_2=25$ (resp. $r_2=21$), then 
$\mu_2=2$ (resp. $\mu_2=4$). Thus $P^2_2=\lbrace (2,25,2),(4,21,2)\rbrace$. 

\par\smallskip

i) $(\mu_2,r_2,d_2)=(2,25,2)$. In this case, applying to construction above to 
$\mu_2=2$, we get 
${\displaystyle{\frac{d_1}{d_2}}=2, {\frac{r_1}{d_2}}=3}$. This leads  to the semigroup $<4,6,25>$. The  
canonical representative of the equisingularity class of this semigroupe is 
$(y^2-x^3)^2-x^{11}y$.

ii) $(\mu_2,r_2,d_2)=(4,21,2)$. In this case, applying to construction above to 
$\mu_2=4$, we get 
${\displaystyle{\frac{d_1}{d_2}}=2, {\frac{r_1}{d_2}}=5}$. This leads 
to the semigroup $<4,10,21>$. The  
canonical representative of the equisingularity class of this semigroupe is 
$(y^2-x^5)^2-x^{8}y$.}

\end{exemple}

\noindent Let $m$ be an even integer, and let $H$ is the set of positive integers not exceeding 
${\displaystyle\frac{{\rm ln}(M)}{{\rm ln}(2)}}$, where 
$M={\displaystyle{\frac{9+\sqrt{1+60m}}{10}}}$. Assume that $H$ is not reduced to $0$
 and let $h$ be a nonzero element of $H$. Set $m=\mu_{h+1}$ and let

$$
a_h={\frac {q+\sqrt{q^2+4\mu_{h+1}.(p+q)}}{2.(p+q)}}=
{\frac {3.2^{h-1}-2+\sqrt{(3.2^{h-1}-2)^2+4\mu_{h+1}.({\frac{5}{3}}.2^{2h-2}-{\frac{2}{3}})}}{{\frac{10}{3}}.2^{2h-2}-{\frac{4}{3}}}}.
$$

\noindent Let $D^h$ be the set  positive integers between 
$2$ and $a_h$ (we easily verify that the condition $a_h\geq 2$ is equivalent to the numerical condition $\mu_{h+1}\geq{\displaystyle  {\frac{5}{3}.2^{2h}-3.2^{h}+{\frac{1}{3}}}}$ proved in Proposition 4.1., in particular $D^h$ is not the emptyset). Set  $P^h=\bigcup_{d\in D^h} P^h_d$ and denote by $C_{h+1}$ the cardinality of $P^h$. In the following we shall give an upper bound for $C_{h+1}$. Set ${\displaystyle b^h_{d}={{\mu_{h+1}\over d}}+
(3.2^{h-1}-1)}$, and  
${\displaystyle c^h_{d}={\mu_{h+1}\over{d-1}}-[{\frac{5}{3}}.2^{2h-2}-3.2^{h-1}+
{\frac{4}{3}}].{d\over {d-1}}+1}$. The set  $R^h_{d}$ of possible values of $r_h$ is a
subset of the set of integers between $b^h_d$ and $c^h_d$.  Its cardinality is then bounded by $c^h_d-b^h_d+1$.Furthermore, we easily verify that
 if $r\in R^h_{d}$, then 

$$
{\displaystyle{\frac{\mu_{h+1}-(r-1)(d-1)}{d}}
\geq {\frac{5}{3}.2^{2h-2}-3.2^{h-1}+{\frac{1}{3}}}}.
$$

\noindent In particular, if $h\geq 2$, then 
$(\mu_h={\displaystyle{\frac{\mu_{h+1}-(r-1)(d-1)}{d}}},r,d)$ is
an element of $P^h_{d}$.

\noindent Now 

$$
c^h_d-b^h_d+1={\frac{\mu_{h+1}}{d(d-1)}}-
(p_h-q_h).{\frac{d}{d-1}}-q_h+1$$

\hskip1cm $$={\frac{\mu_{h+1}}{d-1}}-{\frac{\mu_{h+1}}{d}}-
(p_h-q_h).(1+{\frac{1}{d-1}})-q_h+1
$$

\noindent  Consequently, if  $a=[a_h]$, then the cardinality $C_{h+1}$ of
$P_h$ is bounded by:

$$
{\displaystyle \sum_{d=2}^a
(c^h_d-b^h_d+1)=(\mu_{h+1})(1-{\frac{1}{a}})-p_h(a-1)-(p_h-q_h).\sum_{d=1}^{a-1}
{\frac{1}{d}}+(a-1)}
$$

\noindent But $1-{\displaystyle{\frac{1}{a}} < 1}$, and substituting $2$ to $a$ in the other
members of the above formula we get:  

$$ 
C_{h+1}\leq \mu_{h+1}-2p_h+q_h+1= \mu_{h+1}-({\frac{10}{3}}.2^{2h-2}-3.2^{h-1}-{\frac{1}{3}})=\mu_{h+1}-({\frac{B_{h+1}}{2}}-1)
$$

\noindent Where $B_{h+1}={\displaystyle{\frac{5}{3}}.2^{2h}-3.2^h+{\frac{4}{3}}}$  
is the lower bound of $\mu_{h+1}$ in Proposition 4.1.

\begin{nota}{\rm Note that the bound above is not the optimal one, indeed, given $d\in D^h$,
 the cardinality of the set of $r\in R^h_d$ such that gcd$(r,d)=1$ can be
 bounded by $\displaystyle{\frac{c^h_d-b^h_d}{d}}+1$, but in view of our algorithm,
all values of $R^h_d$ are used, in particular the value above bounds also the
nombre of operations used in the first step of the algorithm.}\end{nota}

\noindent Let $(\mu_h,r,d)$ be an element of $P_h$ and recall that 
$\mu_h={\displaystyle{\frac{\mu_{h+1}}{d}}-(r-1){\frac{d-1}{d}}}$.
Since $b^h_d \leq r \leq c^h_d$, then 

$$
\mu_h \leq {\displaystyle
  {\frac{\mu_{h+1}}{d}}-(b^h_d-1){\frac{d-1}{d}}\leq 
{\frac{\mu_{h+1}}{d}}-(b^h_d-1){\frac{d-1}{d}}={\frac{b^h_d-1}{d}}-(3.2^{h-1}-2)}
$$

\hskip1cm $$={\frac{\mu_{h+1}}{d^2}}+({\frac{1}{d}}-1)(3.2^{h-1}-2)\leq
  {\frac{\mu_{h+1}}{4}}-{\frac{1}{2}}(3.2^{h-1}-2)\leq  {\frac{\mu_{h+1}}{4}}-(3.2^{h-2}-1)
$$

\noindent Let $A_{h+1}=3.2^{h-2}-1$. It follows by induction that for all $0\leq k\leq h-1$, if
$\mu_{h-k}$ is a possible value of the Milnor number at the step $k+1$, then we
have: 

$$
\mu_{h-k}\leq {\frac{\mu_{h+1}}{4^{k+1}}}-\sum_{i=0}^k
{\frac{A_{h+1-i}}{4^{k-i}}}=3.2^{h-2k-2}(2^{k+1}-1)-{\frac{4}{3}}-{\frac{1}{3.4^k}}
$$

\noindent (Remark that the above inequality is valid if $k=h-1$ because $\mu_1=0$). Thus, we can obtain a bound of the set of values 
calculated at the step $k+1, 0\leq k \leq h-2$ in the following way:  let $\mu_{h-k}$ be a possible value of the Milnor
number obtained by reiterating the algorithm above $k+1$ times, and denote by $C_{h-k}(\mu_{h-k})$  the cardinality of the set -denoted $P_{h-k-1}(\mu_{h-k})$- of the 3-uplets $(\mu,r,d)$ obtained by applying the algorithm above to $\mu_{h-k}$ instead of $\mu_{h+1}$. It follows from the discussion above that

$$
C_{h-k}\leq \mu_{h-k}-{\frac{B_{h-k}}{2}}+1\leq {\frac{\mu_{h+1}}{4^{k+1}}}
   -\sum_{i=0}^k {\frac{A_{h+1-i}}{4^{k-i}}}-{\frac{B_{h-k}}{2}}+1$$

$$
={\frac{\mu_{h+1}}{4^{k+1}}}
-3.2^{h-k-2}+3.2^{h-2k-2}-{\frac{5}{3}}.2^{2h-2k-1}-{\frac{1}{3.4^k}}+{\frac{5}{3}}
$$

\noindent In particular, the cardinality of the set of semigroups corresponding to the
given Milnor number $m=\mu_{h+1}$ is bounded by $\prod_{i=2}^{h}C_{h+1-i}$ which is 
a polynomial in $m$ bounded by  its leading coefficient  
$\displaystyle{\frac{m^h}{2^{h(h-1)}}}$. Note that in view of Remark 4.8., the number of
operations used in the algorithm is then bounded by $\sum_{i=0}^{h-1}\prod_{k=0}^{i}C_{h+1-k}$.
 
\medskip

\noindent The above algorithm has been implemented with {\it MAPLE}.
The intput is an integer $m$, and the output is the list of semigroups whose conductor is $m$. 
In the implementation work we followed the ideas explained above, with the following simplification:  at the last step, the set of values we are interested in is calculated by using the factorization of the given Milnor number. The algorithm is an iterating of the following:

{\it Input:  $m\in 2.{\bf N}$

Output: The set $P^h$.

Step I: Compute the set $H$.

Step II:  Take  $h\in H$.

Step III: Compute the set $D^h$.

Step IV: Take $d\in D^h$.

Step V:  Compute $R^h_d$

 \quad\quad (*) if $r\in R^h_d$ and gcd$(r,d)=1$ and $\displaystyle{{\frac{m-(d-1)(r-1)}{d}}\in 2.{\bf N}}$ (resp. $(d-1)(r-1)=m$ if h=1) then add $(\displaystyle{{\frac{m-(d-1)(r-1)}{d}}},r,d)$ to $P^h_d$.

 Step VI: $P^h=\bigcup_{d\in D^h}P^h_d$}

\noindent The main operation of the algorithm  is the one described in the line (*). We 
experimented it on various values of $m$: the computation takes around 0.2sec for $m=160$, 0.7sec for $m=300$, 1.5sec for $m=500$, and 3sec for $m=1000$.

\noindent [1] S.S. Abhyankar.- Lectures on expansion techniques in Algebraic Geometry, Tata Institute
of Fundamental research, Bombay, 1977.

\noindent [2] S.S. Abhyankar.- On the semigroup of a meromorphic curve, Part 1, in Proceedings of
International Symposium on Algebraic Geometry, Kyoto, pp. 240-414, 1977.

\noindent [3] S.S. Abhyankar.- Irreducibility criterion for germs of
analytic functions of two complex variables, Advances in Mathematics
74, pp. 190-257, 1989.

\noindent [4] S.S. Abhyankar.- Some remarks on the Jacobian problem,
Proc.Indian Acad.Sci., vol 104, n$^{0}$ 3, pp. 515-542, 1994.

\noindent [5] S.S. Abhyankar and T.T. Moh.- Newton Puiseux expansion and generalized Tschirnhausen
transformation, J.Reine Angew.Math, 260, pp. 47-83 and 261, pp. 29-54, 1973.

\noindent [6] A. Assi.-Deux remarques sur les racines approch\'ees d'Abhyankar-Moh, C.R.A.S., t.319, Serie 1, 1994, 1191-1196.

\noindent [7] A. Assi.-Meromorphic plane curves, Math.Z., vol 230, 1999, 165-183.

\noindent [8] A. Assi, M. Barile.- Computing irreducible curve singularities with Mathematica, in: L. GonzalézVega, T. Recio (eds.) Actas del 9° Encuentro de Algebra Computacional y
Aplicaciones, Universidad de Cantabria, Santander, 1-3 julio 2004. pp.
17-21.

\noindent [9] A. Campillo.- Algebroid curves in positive characteristic, Lect. Not. Math.,
813. Springer-Verlag, 1980.

\noindent [10] V. Cossart, G. Moreno-Soc{\'\i}as.- Racines approch\'ees 
et suffisance des jets, To appear in Annales de Toulouse.

\noindent [11] V. Cossart, G. Moreno-Soc{\'\i}as.-
Irreducibility criterion: a geometric point of view, in 
Valuation theory and its applications, Fields Institute Communications
Series, AMS, 2002.

\noindent [12] A. Granja.- Irreducible polynomials with coefficients
in a complete discrete valuation field, Advances in Mathematics
109, n 1,  pp. 75-87, 1994. 

\noindent [13] O. Zariski.-Le probl\`eme des modules pour les branches
planes, Hermann, 1986.

\end{document}